\documentclass[12pt,a4paper]{article}
\def\be{\begin{equation}}
\def\ee{\end{equation}}
\def\bea{\begin{eqnarray}}
\def\eea{\end{eqnarray}}
\def\bean{\begin{eqnarray*}}
\def\eean{\end{eqnarray*}}
\def\ben{\begin{enumerate}}
\def\een{\end{enumerate}}
\def\bde{\begin{description}}
\def\ede{\end{description}}
\def\bit{\begin{itemize}}
\def\eit{\end{itemize}}

\begin{document}
\centerline{\textbf{\large  $\lim_{x\to 0}\frac{\sin x}{x}$ and the definition of $\pi$}}
\[
\]
\centerline{\bf Helmut Zeisel}\\
\centerline{e-mail: helmut.zeisel@liwest.at}
\[
\]
\begin{quote}
\begin{abstract}
Leopold Vietoris and Guido Hoheisel showed how the existence of $\lim_{x\to 0}\frac{\sin x}{x}$
can be derived from the trigonometric addition formulas.
In this article two new proofs for this result are given.
In addition it is discussed how this limit is related to the definition of $\pi$.
 
\end{abstract}

{\bf AMS Subject Classification:} 26A09, 39B22, 97D40

{\bf Keywords:} Trigonometric functions, Mathematical induction, Jensen's inequality
\end{quote}

\section{Introduction}
The computation of $\lim_{x\to 0}\frac{\sin x}{x}$ is the fundamental step for the differentiation of the trigonometric functions.
Leopold Vietoris (1957) discussed the usual approaches how to derive this limit and explained the drawbacks of these approaches. As conclusion he showed how the existence of this limit can be derived from the trigonometric addition formulas.
A similar result was found by Hoheisel (1947). These results are also discussed in Acz\'el (1966).

Current textbooks, however, still either use the same old problematic proofs
for $$\lim_{x\to 0}\frac{\sin x}{x}=1$$ or use this equation
as an axiom without further explanation.
For example, Heuser \cite{Heuser} defines sin and cos in axiomatic
way using the addition theorems as functional equations and $\lim_{x\to 0}\frac{\sin x}{x}=1$ as one of the axioms.
He states that the given axioms are not completely independent from each other but does not say in particular how this applies to the limit axiom.

In this article the approach of
Heuser is used. The used axioms to define sine and cosine are, however, slightly different
from Heuers'{}  axioms. Then the
ideas of Vietoris and Hoheisel are applied to show how the existence of $\lim_{x\to 0}\frac{\sin x}{x}$
can be derived from these axioms. 
Finally it is discussed how the value of this limit is related to
the definition of $\pi$.

\section{The definition of the trigonometric functions}
In this section the trigonometric functions are analytically defined by some axioms which can easily be shown geometrically.
These axioms can be considered as a system of functional equations for sine and cosine.

\bea
\cos_c (x \pm y) & = & \cos_c x \cos_c y \mp \sin_c x \sin_c y \\
\sin_c (x \pm y) & = & \sin_c x \cos_c y \pm \cos_c x \sin_c y
\eea
The subscript $c$ is used to point out that it is not specified whether the angles are measured in degrees or radians:
For the moment $c> 0$ is just the angle of the complete circle 
and the right angle has a measure of $c/4$. Additionally the following normalizations are used:
\bea
\cos_c (c/4) & = & 0 \\
\sin_c (c/4) & = & 1
\eea

The particular choice of $c$ is not important
because a change of $c$ is only a change of the angle measurement unit and
can easily be done by the transformations
\be
\sin_{c_1}(x) = \sin_{c_2}\left(\frac{c_2}{c_1}x\right)\mbox{ and }
\cos_{c_1}(x) = \cos_{c_2}\left(\frac{c_2}{c_1}x\right).
\ee

The final axiom is that 
$\sin_c$ is invertible and monotonically increasing and $\cos_c$  is invertible and monotonically decreasing
in the interval $[0, c/4]$. In particular, this implies that both functions are continuous.

If functions satisfying the geometric properties of sine and cosine exist, they must fulfil these axioms.

There are some simple implications from this axioms:
The continuous function $f(x)=\cos_c^2(x)+\sin_c^2(x)$ satisfies the functional equation $f(x+y)=f(x)f(y)$.
This is one of the functional equations discussed by Cauchy \cite{Cauchy}
and the only continuous solution with $f(c/4)=1$ is the function with constant value 1.
\be
\sin_c(0)=\sin_c(x-x)=\sin_c x \cos_c x - \cos_c x \sin_c x=0,
\ee
\be
\cos_c(0)=\cos_c(x-x)=\cos_c^2 x + \sin_c^2 x=1,
\ee
\be
\sin_c(0-x)=\sin_c 0 \cos_c x - \cos_c 0 \sin_c x=-\sin_c x,
\ee
\be
\cos_c(0-x)=\cos_c 0\cos_c x + \sin_c 0\sin_c x=\cos_c x.
\ee

Additionally $\tan_c x$ is defined as $\tan_c x = \frac{\sin_c x}{\cos_c x}$ for $\cos_c x\neq 0$, which is in particular fulfilled for $x\in(-c/4,c/4)$.
Then the law of addition for $\tan_c$ is
\be
\tan_c (x+y) = \frac{\tan_c x + \tan_c y}{1-\tan_c x \tan_c y}.
\ee 
\section{A proof based on induction}
From the above properties, for $x<c/(8n)$ the following inequalities can be derived by induction:
\be
\sin_c (nx) \leq n\sin_c x
\ee
because
\bea
\sin_c ((n+1)x) & = & \sin_c(nx)\cos_c x + \cos_c(nx)\sin_c x \nonumber\\
                & \leq &\sin_c (nx)+ \sin_c x \leq (n+1)\sin_c x. \nonumber
\eea
\be
\tan_c (nx) \geq n\tan_c x
\ee
because
$$\tan_c ((n+1)x) = \frac{\tan_c(nx)+\tan_c x}{1-\tan_c(nx)\tan_c x}\ge n\tan_c x + \tan_c x = (n+1)\tan_c x.$$

Combining these two equations gives
\be
n\sin_c ((n+1)x) \leq (n+1)\sin_c (nx)
\ee
because 
\bean
n \sin_c((n+1)x) & =   & n \tan_c x\cos_c x \cos_c (nx) + n \sin_c (nx) \cos_c x \\ 
               & \le & \tan_c (nx)\cos_c x \cos_c (nx) + n \sin_c (nx) \cos_c x \\
               & \le & (n+1)\sin_c (nx). \nonumber
\eean
and 
\be
n\tan_c ((n+1)x) \geq (n+1)\tan_c (nx)
\ee
because 
\bean
\tan_c((n+1)x) & =   & \frac{n \sin_c (nx)\cos_c x +n \sin_c x \cos_c (nx)}{n\cos_c x \cos_c (nx)-n\sin_c x \sin_c (nx)} \\ 
             & \ge & \frac{n \sin_c (nx)\cos_c x + \sin_c (nx) \cos_c (nx)}{n\cos_c x \cos_c (nx)-\sin_c^2 (nx)} \\
             & =   & \tan_c (nx) \left(1 + \frac{1}{n\cos_c x \cos_c (nx)-\sin_c^2 (nx)}\right)\\
             & \ge & \tan_c (nx) \left(1 + 1/n\right).
\eean

By induction one gets for $m\geq n$
\be
\frac{\sin_c (mx)}{m}\le\frac{\sin_c (nx)}{n}\le\frac{\tan_c (nx)}{n}\le\frac{\tan_c (mx)}{m}
\ee
Now consider two arbitrary positive rational numbers $\frac{p}{q}\leq\frac{r}{s}\leq c/8$,
i.e. $ps\leq qr$ and set $x:=\frac{1}{qs}$, $n:=ps$, and  $m:=qr$.
Then 
\be
\frac{\sin_c (r/s)}{qr}\le\frac{\sin_c (p/q)}{ps}\le\frac{\tan_c (p/q)}{ps}\le\frac{\tan_c (r/s)}{qr}
\ee
and  multiply by $qs$ gives
\be
\frac{\sin_c (r/s)}{r/s}\frac{\le\sin_c (p/q)}{p/q}\le\frac{\tan_c (p/q)}{p/q}\le\frac{\tan_c (r/s)}{r/s}.
\ee
By continuity, for every real $x$, $y$ with $0<x\leq y\leq c/8$
\be\label{monotone}
\frac{\sin_c y}{y}\le\frac{\sin_c x}{x}\le\frac{\tan_c x}{x}\le\frac{\tan_c y}{y}.
\ee

So if $x$ approaches $0$ from the right,
then $\frac{\sin_c x}{x}$ is increasing and bounded from above by the decreasing function $\frac{\tan_c x}{x}$;
in particular $\lim_{x\to 0+}\frac{\sin_c x}{x}$ exists and similarly also $\lim_{x\to 0+}\frac{\tan_c x}{x}$ exists. 
Both limits are nonzero, and, since their quotient is $\cos_c x$, they are equal. Since $\frac{\sin_c x}{x}$ is an even function,  $\lim_{x\to 0+}\frac{\sin_c x}{x}=\lim_{x\to 0-}\frac{\sin_c x}{x}$.

Hoheisel (1947) proved the same inequality (\ref{monotone}) in a different way.

\section{A proof based on convex functions}
Another proof for the differentiability is based on the observation
that $\sin x$ is concave in $[0,c/4]$ and that concave functions have a derivative
everywhere except on at most countable many points, see e.g. \cite{Walter}, p. 304. To prove that
$\sin x$ is concave, it is sufficient to show that it is midpoint concave:

If $f$ is continuous and midpoint convex, i.e. if
$$f\left(\lambda x + (1-\lambda)y\right)\leq \lambda f(x) + (1-\lambda) f(y)$$
holds for $\lambda=\frac{1}{2}$ then it holds for every $\lambda\in[0,1]$, and similarly for concave functions.
This was shown by Jensen \cite{Jensen}. 
A simple proof of Jensen's result can be found in \cite{NiculescuPersson}, Theorem 1.1.4. 

For $0\leq x,y \leq c/8$ 
$$
\sin_c x + \sin_c y = 2 \sin_c \left(\frac{x+y}{2}\right)\cos_c \left(\frac{x-y}{2}\right) \le 2 \sin_c \left(\frac{x+y}{2}\right),
$$
$\sin_c$ is midpoint concave,
\bea
\frac{\tan_c (2x) + \tan_c (2y)}{2}  &-& \tan_c (x+y) \nonumber \\
 & = & \frac{\left(\tan_c x - \tan_c y\right)^2(\tan_c x + \tan_c y)}{(1-\tan_c x \tan_c y)\left(1-\tan_c^2 x\right)\left(1-\tan_c^2 y\right)}\nonumber\geq 0,
\eea
$\tan_c$ is midpoint convex.

Now for $f(0)=0$ and $f$ convex, one has for $0<x\leq y$
$$
\frac{f(x)}{x} = \frac{f(x)-f(0)}{x-0}\le \frac{f(y)-f(0)}{y-0} = \frac{f(y)}{y},
$$
see e.g. \cite{Walter}, p. 303,
and similiar for concave functions, which implies (\ref{monotone}).

\section{\label{pi}The definition of $\pi$}
Up to now it has been shown that $\lim_{x\to 0}\frac{\sin_c x}{x}$ exists;
the value of this limit, however, has not yet been computed.
Hoheisel only mentioned that there exists a value for $c$
such that $\lim_{x\to 0}\frac{\sin_c x}{x}=1$. 
Vietoris used the limit of some ``well known'' recursive sequence for $\pi$
to show  $\lim_{x\to 0}\frac{\sin_{2\pi} x}{x}=1$.
Whether this limit is well known or not depends, however, on the used definition of $\pi$.

The classic geometric definition of $\pi$ is the area of the unit circle, and this area can be computed
as the limit of the areas of inscribed regular $n$-gons. Now the inscribed regular $n$-gon has an area of
$$a_n=\frac{n}{2}\sin_c\frac{c}{n}$$
Substituting $x=c/n$ gives
$$\pi = \lim_{n\to\infty} a_n = \frac{c}{2}\lim_{x\to 0}\frac{\sin_c x}{x},$$
and, as has been shown, this limit exists. Since 
\be
\sin_{c_1}(x) = \sin_{c_2}\left(\frac{c_2}{c_1}x\right)
\ee
this limit is independent of the particular choice of $c$.
So actually there is no need to compute the value of $\lim_{x\to 0}\frac{\sin_c x}{x}$;
it simply can be used as the analytic definition of $\pi$.

Alternatively one could define $\pi$ using the area of the circumscribed
regular $n$-gons, which is
$$A_n=n\tan_c\frac{c}{2n}$$
and leads by defining $x=\frac{c}{2n}$ to 
$$\pi = \lim_{n\to\infty} A_n = \frac{c}{2}\lim_{x\to 0}\frac{\tan_c x}{x},$$
which is the same value as above.

If one chooses $c=2\pi$, one gets $\lim_{x\to 0}\frac{\sin_{2\pi} x}{x}=1,$
and $\sin_{2\pi}$ and $\cos_{2\pi}$ can easily be expanded as a Taylor series.
This corresponds to the often used definition of $\pi/2$ as the first positive zero of $\cos_{2\pi}$. From the presented point of view, however,
this definition of $\pi$ is already implicitly contained in  $\lim_{x\to 0}\frac{\sin_{2\pi} x}{x}=1,$

\section{Summary}
It is possible to prove the exisitence of $\lim_{x\to 0}\frac{\sin x}{x}$
essentially just from trigonometric theorems of addition. This limit can be interpreted as an implicit definition of $\pi$.
This is more intuitive than assuming $\lim_{x\to 0}\frac{\sin x}{x}=1$ just as an axiom
and, in contrast to using geometric proofs, is exact from the point of analysis.

\end{document}